\documentclass[11pt]{article}
\usepackage[dvips]{graphicx}
\usepackage{epsf}
\usepackage{amsfonts}
\usepackage{framed}
\usepackage{setspace}
\newtheorem{theorem}{Theorem}
\newcommand{\nudge}{\hspace{0.01in}}

\begin{document}
\renewcommand{\thefootnote}{\alph{footnote}}
\begin{center}
  {\large \bf CONTINUOUS LUNCHES ARE NOT FREE}\\[.2in] Michael
  D. Vose\nudge\nudge\footnote{EECE Department, University of Tennessee,
    Knoxville, Tennessee 37969, USA.\\ E-mail address:
    vose@eecs.utk.edu}
\end{center}

\date{\pdfcreationdate}

\setcounter{footnote}{0}
\renewcommand{\thefootnote}{\arabic{footnote}}
\begin{quote}
  {\footnotesize A{\tiny BSTRACT}. Stemming from a paper of Auger and
    Teytaud, there is a common misconception that for continuous
    domains No Free Lunch (NFL) does not hold \cite{AT}.  However, Rowe, Vose,
    and Wright have demonstrated that NFL holds for arbitrary domains
    and co-domains \cite{RVW}.  This paper resolves the apparent
    contradiction.  }
\end{quote}


\section{Introduction} The No Free Lunch (NFL) Theorem as introduced
by Wolpert and Macready \cite{WM} originally dealt with Black-box
optimization of functions mapping a finite domain $\mathcal{X}$ to a
finite co-domain $\mathcal{Y}$.  It has since been refined and
generalized to arbitrary $\mathcal{X}$ and $\mathcal{Y}$ \cite{RVW}
and roughly speaking says: {\em All Black-box optimization algorithms
  perform on average equally well over a set $F$ of functions if and
  only if $F$ is closed with respect to permutation}.

In contrast, the paper of Auger and Teytaud \cite{AT} concludes that
NFL does not hold for continuous domains.  As demonstrated below,
their conclusion follows from their imposition of the artificial
constraint -- which is unnecessary for NFL -- that functions under
consideration be measurable.

\section{Measurability and Permutation Closure}

\begin{theorem}
  Let $F$ be a set of Borel measurable functions whose members have
  type $f : [0,1] \rightarrow \mathbb{R}$.  If $F$ is permutation
  closed and $f \in F$, then $f$ is constant almost everywhere.
\end{theorem}
\begin{quote}
  {\em Proof:} Suppose there exists $f \in F$ and $q \in \mathbb{R}$
  such that both inverse images $f^{-1}(-\infty,q]$ and
    $f^{-1}(q,\infty)$ have nonzero Lebesgue measure.  Define $g$ as
  \[
g(x) \; = \; \left\{
\begin{array}{ll}
  1 & \mbox{ if $x > q$}\\
  0 & \mbox{ otherwise}
\end{array}
\right.
\]
and let $h = g \circ f$.  It follows that $h : [0,1] \rightarrow
\{0,1\}$ is Borel measurable, and both preimages $C_0 = h^{-1}\{0\}$
and $C_1 = h^{-1}\{1\}$ are disjoint and uncountable.  Let $S \subset
[0,0.5]$ be a nonmeasurable set; both $S$ and $[0,1] \setminus S$ are
therefore uncountable.  Let $\sigma$ be a bijection $\sigma : S
\rightarrow C_1$ and extend $\sigma$ to a bijection on $[0,1]$ such
that $\sigma : [0,1] \setminus S \rightarrow C_0$.  Note that
\[
h \circ \sigma(x) \; = \; \left\{
\begin{array}{ll}
  1 & \mbox{ if $x \in S$}\\
  0 & \mbox{ if $x \in [0,1] \setminus S$}
\end{array}
\right.
\]
In other words, $h \circ \sigma$ is the indicator function
$\mathbb{I}_S : [0,1] \rightarrow \{0,1\}$ of the nonmeasurable set
$S$.  In particular, $\mathbb{I}_S$ is not Borel measurable.  However,
\[
\mathbb{I}_S \; = \; h \circ \sigma \; = \; g \circ f \circ \sigma \;
= \; g \circ f^\prime
\]
were $f^\prime = f \circ \sigma \in F$ ($F$ is permutation closed).
Because $F$ contains Borel measurable functions, it follows that
$\mathbb{I}_S = g \circ f^\prime$ is Borel measurable, which is a
contradiction.

Therefore, for all $f \in F$ and $q \in \mathbb{R}$, at least one of
\[
L_q = f^{-1}(-\infty,q]\,, \;\;\;\;\; U_q = f^{-1}(q,\infty)
\]
has Lebesgue measure zero. Let $\lambda(\cdot)$ denote Lebesgue
measure. Note that $\lambda(L_q) = \lambda(U_q) = 0$ is impossible,
since $1 = \lambda([0,1]) = \lambda(L_q \cup U_q) = \lambda(L_q) +
\lambda(U_q)$.  Note that
\[
  [0,1] \; = \; f^{-1} \mathbb{R} \; = \; f^{-1} \bigcup_{q \in
    \mathbb{Z}} (-\infty,q] \; = \; \bigcup_{q \in \mathbb{Z}} f^{-1}
    (-\infty,q] \; = \; \bigcup_{q \in \mathbb{Z}} L_q
\]
Hence there exists $q^\prime \in \mathbb{Z}$ such that
$\lambda(L_{q^\prime}) > 0$ and $\lambda(U_{q^\prime}) = 0$.
Similarly,
\[
  [0,1] \; = \; f^{-1} \bigcup_{q \in \mathbb{Z}} (q,\infty) \; = \;
  \bigcup_{q \in \mathbb{Z}} U_q
\]
Hence $\lambda(U_{q^\ast}) > 0$ and $\lambda(L_{q^\ast}) = 0$ for some
$q^\ast \in \mathbb{Z}$.  In particular, $q^\ast < q^\prime$.  Let
$\ell_0 = q^\ast$, $u_0 = q^\prime$ and note that $\lambda(f^{-1}
[\ell_0, u_0]) = 1$.  Define a nested sequence of closed intervals
\[
[\ell_i, u_i] \; \supset \; [\ell_{i+1}, u_{i+1}]
\]
by
\begin{eqnarray*}
m & = & (\ell_{i} + u_{i})/2 \\[.1in]
  [\ell_{i+1}, u_{i+1}] & = & \left\{
\begin{array}{ll}
  [ \ell_{i}, m ] & \mbox{ if $\;\;\lambda(U_m) = 0$} \\[0in]
 [ m, u_{i} ] & \mbox{ if $\;\;\lambda(L_m) = 0$}
\end{array}
\right.
\end{eqnarray*}
It follows that $\lambda(f^{-1} [\ell_i, u_i]) = 1$ for all $i$,
\[
\{c\} \; = \; \bigcap_{i = 0}^\infty \nudge [\ell_i, u_i]
\]
for some $c \in \mathbb{R}$, and $\lambda(f^{-1}\{c\}) = 1$.
In particular, $f$ is constant almost everywhere.
\hfill $\Box$.
\end{quote}

The previous theorem makes clear how the the artificial constraint
that functions be measurable prevents NFL; if $F$ is a collection of
non-trivial (i.e., not constant a.e.) measurable functions, then $F$
cannot be permutation closed!

It could be mentioned that Auger and Teytaud \cite{AT} impose
additional artificial constraints.  However, measurability is enough
to trivialize NFL.

\end{document}